\documentclass{amsart}
\usepackage{amsmath,amssymb}
\usepackage[all]{xy}

\theoremstyle{plain}

\newtheorem{theorem}{Theorem}[section]

\newtheorem{lemma}[theorem]{Lemma}
\newtheorem{corollary}[theorem]{Corollary}

\theoremstyle{definition}

\newtheorem{example}[theorem]{Example}

\theoremstyle{remark}

\newcommand{\secref}[1]{Section~\ref{#1}}
\newcommand{\thmref}[1]{Theorem~\ref{#1}}

\newcommand{\lemref}[1]{Lemma~\ref{#1}}
\newcommand{\corref}[1]{Corollary~\ref{#1}}

\newcommand{\exref}[1]{Example~\ref{#1}}

\def\R{{\mathbb \R}}
\def\Z{{\mathbb Z}}
\def\Q{{\mathbb Q}}
\def\C{{\mathbb C}}

\def\map{\mathrm{map}}

\def\Hom{\mathrm{Hom}}

\def\cat0{\mathrm{cat}_0}

\begin{document}

\title[Components of Function Spaces]
{Criteria for Components of a Function Space to be Homotopy
Equivalent}

\author{Gregory  Lupton}

\address{Department of Mathematics,
	   Cleveland State University,
	   Cleveland OH 44115}

\email{G.Lupton@csuohio.edu}

\author{Samuel Bruce Smith}

\address{Department of Mathematics,
   Saint Joseph's University,
   Philadelphia, PA 19131}

\email{smith@sju.edu}

\date{\today}



\begin{abstract}
We give a general method that may be effectively applied to the
question of whether two components of a function space $\map(X,Y)$
have the same homotopy type.  We describe certain group-like actions
on $\map(X,Y)$. Our basic results assert that if maps $f, g \colon X
\to Y$ are in the same orbit under such an action, then the
components of $\map(X, Y)$ that contain $f$ and $g$ have the same
homotopy type.
\end{abstract}

\maketitle

\section{Introduction}%
\label{sec:intro} Let $X$ and $Y$ be connected, countable CW
complexes and let $\map(X, Y)$ denote the space of all continuous
(not necessarily based) maps between $X$ and $Y$ with the
compact-open topology.  The space $\map(X, Y)$ is generally
disconnected with path components in one-to-one correspondence with
the set $\langle X,Y\rangle$ of (free) homotopy classes of maps.
Furthermore, different components may---and frequently do---have
distinct homotopy types.  A basic problem in homotopy theory is to
determine whether two components are homotopy equivalent or, more
generally, to classify the path components of $\map(X, Y)$ up to
homotopy equivalence. For $x_0 \in X$ a choice of basepoint, we have
the evaluation map $\omega \colon \map(X,Y) \to Y$, defined by
$\omega(g) = g(x_0)$, which is a fibration.  Let $\map(X, Y;f)$
denote the path component of $\map(X, Y)$ that contains a given map
$f \colon X \to Y$.  We may also ask for a finer classification, up
to fibre-homotopy equivalence, of the evaluation fibrations
$\omega_f \colon \map(X, Y;f) \to Y$, obtained by restricting
$\omega$ to the component of $f$.

Work on these classification problems dates back to the 1940s.
Whitehead considered the case $X=S^n$ and $Y=S^m$, in which a
component corresponds to $\alpha \in \pi_n(S^m)$, and proved that
$\map(S^n, S^m; \alpha)$ is homotopy equivalent to $\map(S^n, S^m;
0)$ if and only if the evaluation fibration $\omega_\alpha$ admits a
section \cite[Th.2.8]{Wh1}. Hansen, and later McLendon, extended
this analysis (\cite{Han1, Han2, McCl}). In \cite{Han3}, Hansen
obtained a classification of components of $\map(M^n,S^n)$, where
$M^n$ is a suitably restricted $n$-manifold.  Sutherland extended
this result in \cite{Suth1}. M{\o}ller \cite{Moller} gave a
classification of components of $\map(\C P^m, \C P^n)$ for $1 \leq m
\leq n$.  The case in which $X$ is a manifold and $Y = BG$, the
classifying space of a Lie group, has been the subject of extensive
recent research by Crabb, Kono, Sutherland, Tsukuda and others (see
e.g. \cite{ Crabb-Suth,Kono, K-T, K-T2,Suth2,Ts}). Our purpose in
this paper is to give a general method that may be applied to show
that (evaluation fibrations of) components of $\map(X, Y)$ are
(fibre-) homotopy equivalent. In addition to yielding many new
results, our method allows some of the particular cases just
mentioned to be viewed as special cases within a general framework.

Our basic results are presented in \secref{sec:action}.  We consider
the orbit of a point in $\map(X,Y)$ under a group-like action on
$\map(X,Y)$ and observe in \thmref{thm:components of an action} that
two distinct components of $\map(X,Y)$ are homotopy equivalent
whenever each overlaps with any one orbit---not in the same point,
obviously. Now, in the situations that we have in mind, the action
on $\map(X,Y)$ arises from a group-like action on $Y$. In this case,
we have a corresponding group action on the set of homotopy classes
of maps $\langle X,Y\rangle $. Write $\mathcal{O}$ for the orbit set
of this group action. Then we obtain a surjection
\begin{equation}\label{eq:surjection}
\xymatrix{ \mathcal{O} \ar@{->>}[rr] &&
\begin{displaystyle}
\frac{\{ \text{components of $\map(X,Y)$} \}}{\simeq}
\end{displaystyle}}
\end{equation}
of sets, where $\simeq$ denotes homotopy equivalence
(\thmref{thm:actions on Y}).  This may be applied ``locally," to
analyze whether two particular components are homotopy equivalent.
It may also be applied ``globally," to deduce a finite---or even a
concrete upper bound on the---number of distinct homotopy types
amongst the (usually infinitely many) components of $\map(X,Y)$. We
illustrate both approaches in \secref{sec:holonomy}.  For based
spaces $X$ and $Y$, we may also consider $[X,Y]$, the set of
based-homotopy equivalence classes of based maps. Ignoring
basepoints gives a surjection $\xymatrix{ [X,Y] \ar@{->>}[r]
&\langle X,Y\rangle  }$ of sets of homotopy classes. Once more, in
the situations that we have in mind, the group action on $\langle
X,Y\rangle $ that we referred to above actually restricts to one on
$[X,Y]$. Writing $\mathcal{O}_*$ for the corresponding orbit set, we
may compose the surjection (\ref{eq:surjection}) with the surjection
$\xymatrix{ \mathcal{O}_* \ar@{->>}[r] &\mathcal{O} }$ of orbit
sets.  Although $\mathcal{O}_*$ is \emph{a priori} larger than
$\mathcal{O}$, it is more familiar in homotopy theory and in many
cases may be analyzed effectively. With further restrictions on $X$
and $Y$, we may sharpen these results, replacing the right-hand set
in (\ref{eq:surjection}) by fibre-homotopy equivalence classes of
evaluation fibrations $\omega_f \colon \map(X, Y;f) \to Y$. Also, we
may readily adapt the methods used here to study homotopy types of
components of $\map_*(X, Y)$, the function space of
basepoint-preserving maps---see the comment at the end of
\secref{sec:action} and the discussion that ends the paper.

In \secref{sec:holonomy} we focus our general method on actions on
$\map(X,Y)$ that arise from certain specific actions on $Y$. We
first consider the holonomy action of $\Omega B$ on the fibre $Y$ of
a fibration $Y \to E \to B$. In \thmref{thm:pre-image of j_*}, we
show that if two based maps $f, g \colon X \to Y$ satisfy $j\circ
f\sim_* j\circ g \colon X \to E$, where $j \colon Y \to E$ denotes
the fibre inclusion, then the components $\map(X,Y;f)$ and
$\map(X,Y;g)$ have the same homotopy type. With some restrictions on
$X$ and $Y$, we are able to conclude more strongly that the
evaluation fibrations $\omega_f \colon \map(X, Y; f) \to Y$ and
$\omega_g \colon \map(X, Y;g) \to Y$ are fibre-homotopy equivalent.
We illustrate these ideas in \exref{ex:G/T} and \exref{ex:BH}, which
give simple cohomological conditions under which two components of
$\map(X,G/H)$ are homotopy equivalent, or there are finitely many
homotopy types amongst the components of $\map(X,G/H)$, where $H$ is
a closed subgroup of a Lie group $G$. Next we focus on the universal
fibration with fibre $Y$ and obtain a link between the
classification problem for components of a function space and the
class of \emph{cyclic maps} (see \cite{Var}). In this context, we
extend the result of Whitehead mentioned above to prove that the
evaluation fibrations $\omega_f \colon \map(X, Y; f) \to Y$ and
$\omega_0 \colon \map(X, Y;0) \to Y$ are fibre-homotopy equivalent
if and only if $\omega_f$ admits a section (\thmref{thm:cyclic}). We
obtain further results in the case in which $X$ is a co-H-space,
including a connection between computations of the Gottlieb groups
of spheres and Hansen's results on the classification of the
components of $\map(S^n, S^m)$ (cf.~\exref{ex:Stiefel manifold}). We
end the paper with a brief discussion of comparable results about
components of the based mapping space $\map_*(X, Y)$, but with the
action arising from cogroup-like actions on $X$.

\section{Group-Like Actions on a Function Space}%
\label{sec:action}

We begin by setting conventions and notation. First, we make clear
that \emph{homotopy} (homotopic maps, homotopy equivalence, etc.)
generally refers to \emph{free homotopy}: we use ``$\sim$" and
``$\simeq$" to denote (free) homotopy and (free) homotopy
equivalence, respectively.  If \emph{based homotopy} is intended, we
will be specific and use ``$\sim_*$" and ``$\simeq_*$" in that case.

A \emph{fibration} $p \colon E \to B$ means a Hurewicz fibration
\cite[p.29]{Wh}. Recall that, for $p_1 \colon E_1 \to B$ and $p_2
\colon E_2 \to B$ fibrations over a space $B$, a based map $f \colon
E_1 \to E_2$ is a \emph{fibre homotopy equivalence} if there exists
$g \colon E_2 \to E_1$ such that $g \circ f $ and $f \circ g$ are
homotopic to the respective identities by based homotopies $F$ and
$G$ satisfying $p_1\circ F(x, t) = p_1(x)$ and $p_2\circ G(y, t) =
p_2(y)$ for $x \in E_1, y \in E_2$ and $t \in I.$

An \emph{H-space} is a based space $G$ together with a based
multiplication $m \colon G \times G \to G$ that satisfies $m\circ J
\sim_*  \nabla\colon G \vee G \to G$ where $J \colon G \vee G \to G
\times G$ is the inclusion and $\nabla \colon G \vee G \to G $ is
the folding map.   We note that the homotopy can be replaced by
strict equality provided the basepoint of $G$ is non-degenerate
\cite[Thm.III.4.7]{Wh}. The H-space is \emph{homotopy-associative}
if $m\circ(m\times 1)\sim_* m\circ(1\times m)\colon G\times G\times
G\to G$. By a \emph{group-like} space, we mean a
homotopy-associative H-space $G$ together with a based inverse map
$\iota \colon G \to G$ that satisfies $m\circ(\iota\times 1)\circ
\Delta \sim_* 0$ and $m\circ(1\times\iota)\circ \Delta \sim_* 0$,
where $\Delta \colon G \to G \times G$ is the diagonal map.

By a \emph{homotopy-associative action} of a homotopy-associative
H-space $G$ on a based space $Y,$ we mean a based map $A \colon G
\times Y \to Y$ that satisfies $A\circ i_2 \sim_* 1\colon Y \to Y$
and $A\circ(1\times A) \sim_* A\circ (m\times 1)\colon G\times G
\times Y \to Y$, where $i_2 \colon G \to Y \times G$ is the
inclusion. We say the action is \emph{strictly unital} if we have $A
\circ i_2 = 1.$ The argument in \cite[Thm.III.4.7]{Wh} mentioned
above easily extends to show an action may be taken to be strictly
unital when the basepoint of $G$ is non-degenerate. Given $g \in G$
and $x \in Y$, we will usually write $g\cdot x$ for $A(g, x)$.

For the rest of the paper, we assume (at least) that $X$ and $Y$ are
based, connected, countable CW complexes with fixed choices of
non-degenerate basepoints.  While these hypothesis are not strictly
necessary for all that we do, they seem to provide a reasonable
level of generality.  Despite these restrictions on $X$ and
$Y$---indeed, despite further restrictions (e.g. $X$ is frequently
assumed to be a finite complex)---we must allow for much greater
generality when considering the function space $\map(X,Y)$.
\lemref{lemma:nondegenerate} and \lemref{lem:non-deg map(X,Y)} below
deal with technical points that become issues when we consider the
function space.

\begin{lemma}\label{lemma:nondegenerate}
Suppose $U$ and $V$ are path-connected spaces with non-degenerate
basepoints.  Then we have:
\begin{itemize}
\item[(1)] Given $f \colon U \to V$, there
exists a based map $f' \colon U \to V$ with $f \sim f'$;
\item[(2)] If $U$ and $V$ are homotopy equivalent, then they are
based homotopy equivalent with respect to the non-degenerate
basepoints.
\end{itemize}
\end{lemma}

\begin{proof} Part (1) is \cite[III.1.4]{Wh}. For (2),
suppose that $f \colon U \to V$ is a (free) homotopy equivalence.
Let $u_0 \in U$ and $v_0 \in V$ be non-degenerate basepoints. Since
$u_0$ is non-degenerate, $f$ is homotopic to a based map $f'\colon U
\to V$ by (1). Since $f$ is a homotopy equivalence, so too is $f'$.
But since $f'(u_0) = v_0$, and both $u_0$ and $v_0$ are
non-degenerate, it follows that $f'$ is a based homotopy equivalence
(see, e.g., \cite[Prop.6.18]{Jam84}).
\end{proof}

It is well-known that all components of a group-like space have
the same homotopy type (see \cite[Prop.5.28]{Jam84}). We
generalize this fact in the following result:

\begin{theorem}\label{thm:components of an action}
Let $A \colon G \times Y \to Y$ be a homotopy-associative action of
a group-like space $G$ on a space $Y$.  For each $x \in Y$, let $Y_x
\subseteq Y$ denote the path component of $Y$ that contains $x$.
Then for each $g \in G$, the components $Y_x$ and $Y_{g\cdot x}$
have the same homotopy type.  If $Y_x$ and $Y_{g\cdot x}$ both have
non-degenerate basepoints, then $Y_x$ and $Y_{g\cdot x}$ have the
same based homotopy type.
\end{theorem}

\begin{proof}
Let $m \colon G\times G \to G$ be the multiplication and
$\iota\colon G \to G$ the inverse map.  Let $e \in G$ denote the
basepoint. For each $g \in G$, we may define ``translation by $g$"
to be the map $\tau_g \colon Y \to Y$, where $\tau_g(x) = g\cdot
x$ for each $x \in Y$. Then $\tau_g$ restricts to a map $\tau_g
\colon Y_x \to Y_{g\cdot x}$. On the other hand, we have the
translation $\tau_{\iota(g)}\colon Y \to Y$. Let $i_g \colon Y \to
G\times Y$ be the inclusion defined by $i_g(x) = (g, x)$ for each
$x \in Y$. Then we have
\begin{equation*}
\begin{aligned}
\tau_{\iota(g)}\circ \tau_{g} &= A\circ (1\times A) \circ\big(
(\iota, 1) \times 1\big)\circ i_g \sim A\circ (m\times 1)
\circ\big( (\iota, 1) \times 1\big)\circ i_g\\
&\sim A\circ (0 \times 1)\circ i_g \sim \tau_e \sim 1_Y.
\end{aligned}
\end{equation*}
Let $H \colon Y \times I \to Y$ be a homotopy from
$\tau_{\iota(g)}\circ \tau_{g}$ to $1$. $H(x, t)$ gives a path from
$\iota(g)\cdot(g\cdot x)$ to $x$ and it follows that
$\tau_{\iota(g)}$ restricts to a map $\tau_{\iota(g)}\colon
Y_{g\cdot x} \to Y_x$.  Furthermore, the homotopy $H$ restricts to a
homotopy $H \colon Y_x\times I \to Y_x$ between the composition of
the restrictions $\tau_{\iota(g)}\circ \tau_{g}$ and the restriction
of the identity to $Y_x$.  That is, the restriction of
$\tau_{\iota(g)}$ to $Y_{g\cdot x}$ is a left-homotopy inverse for
the restriction of $\tau_g$ to $Y_x$.  A similar argument shows that
$\tau_{\iota(g)}$ is a two-sided inverse, and thus $\tau_g \colon
Y_x \to Y_{g\cdot x}$ is a homotopy equivalence. The last assertion
follows from \lemref{lemma:nondegenerate}~(2).
\end{proof}

Our interest in \thmref{thm:components of an action} lies in its
implications for function spaces.  By
\lemref{lemma:nondegenerate}~(1), any map $f \colon X \to Y$ is
homotopic to a based map. Therefore, when identifying a component of
$\map(X,Y)$ as $\map(X,Y;f)$ for some map $f \colon X \to Y$, we may
assume that $f$ is a based map. Also, these hypotheses ensure that
the evaluation map $\omega_f \colon \map(X, Y;f) \to Y$ is a
Hurewicz fibration by \cite[Th.I.7.1]{Wh}. We write $\map_*^*(X,
Y;f) =\omega_f^{-1}(*)$ for the fibre over the basepoint of $Y$.
Note that the space $\map_*^*(X, Y;f)$ consists of based maps $g
\colon X \to Y$ which are (freely) homotopic to $f$. Thus $\map_*(X,
Y;f) \subseteq \map_*^*(X, Y;f)$ and the inclusion can be strict.

\lemref{lemma:nondegenerate} indicates that we will want $\map(X,Y)$
to have non-degenerate basepoints.  Since we have not been able to
find an explicit reference for what we want in the literature, we
provide the following result that is suited to our purposes.

\begin{lemma}\label{lem:non-deg map(X,Y)}
Let $X$ be a compact metric space and $Y$ a countable CW complex.
Then every point of $\map(X,Y)$ is non-degenerate.
\end{lemma}

\begin{proof}
Suppose given $a_0 \in \map(X,Y)$.  We want to show that $\{ a_0 \} 
\hookrightarrow
\map(X,Y)$ is a cofibration.  By \cite[Th.XV.7.4]{Dugundji}, it is
sufficient to show that there exists a neighbourhood of $a_0$ in
$\map(X,Y)$, of which $\{ a_0 \}$ is a strong deformation retract.  By
\cite[Lem.3]{Mil}, $\map(X,Y)$ is ``ELCX."  In particular, we may
choose an open set $V$ in $\map(X,Y)$ that contains $a_0$, and for
which there exists a homotopy $\lambda \colon V \times V \times I
\to \map(X,Y)$ that satisfies $\lambda(a, b, 0) = a$, $\lambda(a, b,
1) = b$, and $\lambda(a, a, t) = a$ for all $a, b \in V$ (cf.~the
discussion above \cite[Lem.3]{Mil}). We define $H \colon V \times I
\to \map(X,Y)$ by $H(b,t) = \lambda(a_0, b, t)$ and check that this
displays $\{a_0 \}$ as a strong deformation retract of $V$, as required.
\end{proof}

Suppose $A \colon G \times Y \to Y$ is a homotopy-associative
action of a group-like space $G$   on $Y$.
The based function space $\map_*(X, G)$ is then a group-like space
as well,  with multiplication defined to be  pointwise multiplication of
functions \cite[Th.III.5.18]{Wh}.  We have an induced action
$$\mathcal{A} \colon \map_*(X, G) \times \map(X, Y) \to \map(X, Y)$$
of $\map_*(X, G)$ on $\map(X, Y)$ defined by
$$\mathcal{A}(\gamma,
g)(x) = A(\gamma(x), g(x))$$
for $\gamma \in \map_*(X, G), g \in \map(X, Y).$ As above, we write
$\gamma \cdot g$  for $\mathcal{A}(\gamma, g)$.  We note that the
following result holds in considerable generality.
\begin{theorem}\label{thm:actions on Y}
Let $f \colon X \to Y$ be a map between based, connected, countable
CW complexes. Let $A \colon G \times Y \to Y$ be  a
homotopy-associative action of a group-like space $G$  on $Y$.  Let
$\gamma \colon X \to G$ be any based map. Then we have:
\begin{itemize}
\item[(A)] The path components $\map(X,Y;f)$ and
$\map(X,Y;\gamma\cdot f)$ have the same homotopy type;
\item[(B)] If $X$ is a finite complex, then
$\map(X,Y;f)$ and $\map(X,Y;\gamma\cdot f)$ have the same based
homotopy type.
\item[(C)] If $X$ is finite and the action is strictly unital,
then the evaluation fibrations $\omega_f \colon \map(X, Y;f) \to Y$
and $\omega_{\gamma \cdot f} \colon \map_{\gamma \cdot f}(X, Y;
\gamma \cdot f)$ are fibre-homotopy equivalent.
\end{itemize}
\end{theorem}
\begin{proof}
Write $\tau_\gamma \colon \map(X, Y;f) \to \map(X, Y; \gamma \cdot
f)$ for translation by $\gamma.$  By \thmref{thm:components of an
action}, $\tau_\gamma$ is a homotopy equivalence and (A) follows.
Part (B) follows from \thmref{thm:components of an action} and
\lemref{lem:non-deg map(X,Y)}. For (C), we use results of Dold in
\cite{Dold}. The evaluation fibration $\omega_f \colon \map(X, Y;f)
\to Y$ has the Weak Covering Homotopy Property \cite[Def.5.1]{Dold}
since it is a Hurewicz fibration. Since the action of $G$ is
strictly unital, the diagram
$$ \xymatrix{  \map(X, Y;f)  \ar[rd]_{\omega_f}
\ar[rr]^{\tau_\gamma}  && \map(X, Y;\gamma \cdot f)
\ar[dl]^{\omega_{\gamma \cdot f}} \\
& Y }$$ commutes. The map $\tau_\gamma$ is a based homotopy
equivalence by (B). By \cite[Th.6.1]{Dold}, $\tau_\gamma$ is thus a
fibre-homotopy equivalence.
\end{proof}

We can recast \thmref{thm:actions on Y} as follows: Write
$$\mathcal{A}_\# \colon [X, G] \times \langle X,Y\rangle \to \langle X,Y\rangle$$
for the induced action of the group $[X, G]$ induced on the set
$\langle X,Y\rangle  $ of homotopy classes of maps. We write
$\mathcal{O}$ for the set of orbits of $\langle X,Y\rangle  $ under
this action.
\begin{corollary} \label{cor:surjection}
Let  $f \colon X \to Y$ be a map between CW complexes. Let $A \colon
G \times Y \to Y$ be  a homotopy-associative action of a group-like
space $G$  on $Y$. Let $\mathcal{O}$ be the set of orbits of the
induced action of the group $[X, G]$ on $\langle X,Y\rangle .$ Then
\begin{itemize}
\item[(A)] We have a surjection of sets
$$ \xymatrix{ \mathcal{O} \ar@{->>}[rr] &&
\begin{displaystyle}
\frac{\{ \text{components of $\map(X,Y)$} \}}{\simeq}
\end{displaystyle}.}$$
\item[(B)]
If $X$ is a finite complex then
$$ \xymatrix{ \mathcal{O} \ar@{->>}[rr] &&
\begin{displaystyle}
\frac{\{ \text{components of $\map(X,Y)$} \}}{\simeq_*}
\end{displaystyle}.}$$
\item[(C)]
If $X$ is a finite complex and the group-like action  on $Y$ is strictly unital then
$$ \xymatrix{ \mathcal{O} \ar@{->>}[rr] &&
\begin{displaystyle}
\frac{\{ \text{evaluation fibrations   $\omega_f \colon \map(X,Y;f)
\to Y$} \}}{\text{fibre-homotopy equivalence}}
\end{displaystyle}.}\qed$$
\end{itemize}
\end{corollary}

In particular, if $\mathcal{O}$ is a finite set, then there are
finitely many distinct homotopy types amongst the components of
$\map(X,Y)$.

We observe that the discussion of this section can be given with
$\map_*(X,Y)$ replacing $\map(X,Y)$.  We will see in the next
section that there is a further situation that gives rise to an
action on $\map_*(X,Y)$, to which we may apply our methods.

\section{Holonomy Actions and Universal Actions}\label{sec:holonomy}

A standard source for an action on a space $Y$ is fibration sequence
\begin{equation}\label{eq:general fibration}
\xymatrix{Y \ar[r]^{j} & E \ar[r]^{p} & B}
\end{equation}
in which $Y$ occurs as the fibre.  For then we have the holonomy
action $A\colon \Omega B \times Y \to Y$ of the group-like space
$\Omega B$ on $Y$.  As above, this yields an induced action
\begin{equation}\label{eq:action induced from fibration}
\mathcal{A} \colon \map_*(X, \Omega B) \times \map(X, Y) \to
\map(X, Y)
\end{equation}
of $\map_*(X, \Omega B)$ on $\map(X, Y)$.   In this situation, we
may be quite precise about the orbits.

\begin{lemma}
\label{lemma:holonomy}
Let $f, g \colon X \to Y$ be based maps.  With reference to the
action (\ref{eq:action induced from fibration}) induced from the
fibration (\ref{eq:general fibration}), the following are
equivalent:
\begin{itemize}
\item[(A)] $g \sim_* \gamma\cdot f$ for some
$\gamma \in \map_*(X, \Omega B)$, that is,
$f$ and $g$ are in the same orbit;
\item[(B)] $g \sim_* A\circ(\gamma\times f)\circ\Delta$;
\item[(C)] $j\circ f \sim_* j\circ g \colon X \to E$.
\end{itemize}
\end{lemma}

\begin{proof}
(A) and (B) are equivalent from the definitions. To see (C) is
equivalent, consider the Puppe sequence
\begin{equation}\label{eq:Puppe sequence fibration}
\xymatrix{ \cdots \ar[r] & [X, \Omega B] \ar[r]^{\partial_*} & [X,
Y] \ar[r]^{j_*} & [X, E] \ar[r]^{p_*} & [X, B] }
\end{equation}
corresponding to the fibration (\ref{eq:general fibration}).  As
is well-known, $[X, \Omega B]$ acts on $[X, Y]$ as described in
(B) (see e.g. \cite[p.140]{Wh}).  Furthermore, an orbit of $[f] \in [X,Y]$ under this action
is precisely the pre-image of $j_*([f])$ (\cite[III.6.20]{Wh}).  The equivalence of (B)
and (C) follows.
\end{proof}

\begin{theorem}\label{thm:pre-image of j_*}
Let $X$ and $Y$ be connected, countable CW complexes with
non-degenerate basepoints. Let $j \colon Y \to E$ be the fibre
inclusion of a fibration in which $Y$ occurs as the fibre.  Suppose
$j\circ  f \sim_* j\circ g \colon X \to E$ for maps $f, g \colon X
\to Y$. Then $\map(X,Y;f)$ and $\map(X, Y; g)$ are homotopy
equivalent. If $X$ is a finite complex then the evaluation
fibrations $\omega_f \colon \map(X, Y;f) \to Y$ and $\omega_g \colon
\map(X, Y;g) \to Y$ are fibre-homotopy equivalent.
\end{theorem}

\begin{proof}
The result follows directly from \thmref{thm:actions on Y} and
\lemref{lemma:holonomy}.
\end{proof}

\begin{example}\label{ex:G/T}
Consider a compact, connected Lie group $G$ and a toral subgroup $T
\subseteq G$. Then we have a fibre sequence $G/T \to BT \to BG$ with
fibre inclusion $j\colon G/T \to BT = \prod K(\Z,2)$. Given a CW
complex $X$ and based maps $f, g \colon X \to G/T$ we see that
$H^2(f) = H^2(g) \colon H^*(G/T;\Z) \to H^*(X; \Z)$ implies $j \circ
f \sim_* j \circ g$.  We conclude from \thmref{thm:pre-image of j_*}
that $H^2(f) = H^2(g)$ implies the components $\map(X, G/T; f)$ and
$\map(X, G/T; g)$ are homotopy equivalent.
\end{example}

We may develop \thmref{thm:pre-image of j_*} as follows.

\begin{corollary}\label{cor:j_*}
Let $X$ and $Y$ be connected, countable CW complexes with
non-degenerate basepoints. Let $j \colon Y \to E$ be the fibre
inclusion of a fibration in which $Y$ occurs as the fibre. If the
image of $j_* \colon [X,Y] \to [X, E]$ is a finite set in $[X, E]$,
then there are finitely many distinct homotopy types amongst the
components of $\map(X,Y)$.  If, further, $X$ is finite, there are
finitely many fibre-homotopy types amongst the evaluation fibrations
$\omega_f \colon \map(X, Y;f) \to Y$ for $f \colon X \to Y.$ \qed
\end{corollary}

So, for instance, returning to the situation of \exref{ex:G/T}, we
may say that \emph{if $H^2(X; \Z)$ is finite, then there are
finitely many distinct homotopy types amongst the components of
$\map(X,G/T)$.} We offer a further example along these lines.

\begin{example}\label{ex:BH}
Let $G$ be a connected Lie group and $H$ a closed subgroup.  Suppose
that $\Hom\big( H^*(G/H;\Q), H^*(X; \Q)\big) = 0$, for a finite
complex $X$.  (These hypotheses hold, for instance, whenever $H$ is
a subgroup of maximal rank and $X$ is any finite complex with
$H^{\mathrm{even}}(X;\Q) = 0$.) Then there are finitely many
fibre-homotopy types amongst the evaluation fibrations $\omega_f
\colon \map(X, G/H;f) \to Y$, for maps $f \colon X \to G/H$.  To see
why, observe that $BH$ is rationally a product of
Eilenberg-Mac\;Lane spaces, and hence the hypotheses imply that each
$j\circ f \colon X \to BH$ is null-homotopic after rationalization,
where $j\colon G/H \to BH$ is the fibre inclusion of the fibre
sequence $G/H \to BH \to BG$ and $f \colon X \to G/H$ is any map.
Since rationalization of homotopy sets is a finite-to-one map
\cite[Cor.II.5.4]{H-M-R}, it follows that $j_* \colon [X, G/H] \to
[X, BH]$ has finite image.  Now we may apply \corref{cor:j_*}.
\end{example}

We next observe that the universal action on a space $Y$ is that
induced by the evaluation map of the identity component. Precisely,
observe that the space $\map(X, X; 1)$ is a strictly associative
$H$-space with multiplication given by composition of functions.
Define the action $$A_\infty \colon \map(Y, Y;1) \times Y \to Y$$ by
$A_\infty(g, y) = g(y)$ for $g \in \map(Y, Y;1)$ and $y \in Y.$
Given any $H$-action $A \colon G \times Y \to Y$ we obtain, by
adjointness, an $H$-map $\hat{A} \colon G \to \map(Y, Y; 1)$ which
commutes with the actions in the sense that $A_\infty(\hat{A}(g), y)
= A(g, y)$ for all $g \in G$ and $y \in Y$. Conversely, any $H$-map
$\hat{A} \colon G \to \map(Y, Y; 1)$ induces an action $A \colon G
\times Y \to Y$.  We remark that, according to Gottlieb \cite{Got2},
this universal action corresponds to the holonomy action in the
universal fibration with fibre $Y$
\begin{equation}\label{eq:classifying fibration}
\xymatrix{Y \ar[r]^-{j_\infty} & E_\infty \ar[r]^-{p_\infty} &
B_\infty,}
\end{equation}
(cf.~\cite{St, All, Dold1}).  We will need the following consequence
of the classifying fibration:
\begin{theorem} Let $Y$ be a CW complex.  Then $\map(Y, Y;1)$ is a
group-like space.
\end{theorem}
\begin{proof}
Since $Y$ is a CW complex, combining \cite{All} with
\cite[Satz.7.3]{Fuchs} gives an $H$-equivalence $\map(Y, Y;1) \simeq
\Omega_0 B_\infty$, where $\Omega_0 B_\infty$ denotes the component
of the constant loop.  Thus $\map(Y, Y;1)$ is group-like by
\cite[Cor.III.5.17]{Wh}.
\end{proof}

Now write
$$\mathcal{A}_\infty \colon \map_*(X, \map(Y, Y;1)) \times \map(X, Y)
\to \map(X, Y)$$
for the action induced by $A_\infty$ on $\map(X,Y)$, and
$$(\mathcal{A}_\infty)_\sharp \colon [X,   \map(Y,Y;1)]
\times [X, Y] \to [X, Y]$$
for the corresponding group action on the set $[X, Y].$   Given a
homotopy class $[f] \in [X, Y]$ we write $O_\infty([f])$ for the
orbit of $[f]$ under this action.

We recall that a based map $f \colon X \to Y$ is called
\emph{cyclic} if the map $(f\mid 1) \colon X \vee Y \to Y$ admits
some extension $\Gamma \colon X \times Y \to Y$ \cite{Var}. We write
$G(X, Y) \subseteq [X, Y]$ for the set of based homotopy classes of
cyclic maps. In the special case in which $X =S^n$, $G(S^n, Y)$ is
just $G_n(Y) \subseteq \pi_n(Y)$, the \emph{$n$th Gottlieb group of
$Y$} \cite{Got3}.

It is a direct consequence of adjointness that $f \colon X \to Y$ is
cyclic if and only if the evaluation fibration $\omega_f \colon
\map(X, Y;f) \to Y$ admits a section. We also have
\begin{equation} \label{eq:orbit=cyclic}
 O_\infty([0]) = G(X,Y).
\end{equation}
For suppose $\gamma \colon X \to \map(Y,Y;1)$ is a based map. We
define a section $s \colon Y \to \map(X, Y;\gamma \cdot 0)$ by the
rule $s(y) = \gamma(x)(y).$ Conversely, if $s \colon Y \to \map(X,
Y;f)$ is a section for some based map $f \colon X \to Y$ then $f
\sim_* \gamma \cdot 0$ where $\gamma \colon X \to \map(Y, Y;1)$ is
given by $\gamma(x)(y) = s(y)(x).$ As a consequence,   we obtain the
following result which extends \cite[Th.2.8]{Wh1} and its
generalization by Yoon \cite[Th.4.5]{Yoon89}.

\begin{theorem}
\label{thm:cyclic} Let $X$ and $Y$ be CW complexes with
non-degenerate basepoints.  Let $f \colon X \to Y$ be a map. If
$\omega_f \colon \map(X, Y;f) \to Y$ has a section then $\map(X,
Y;f)$ is homotopy equivalent to $\map(X, Y;0)$.   If $X$ is a finite
complex then the following are equivalent:
\begin{itemize}
\item[(A)] The map $f \colon X \to Y $ is cyclic.
\item[(B)] The evaluation fibration $\omega_f \colon \map(X, Y;f)
\to Y$
has a section.
\item[(C)] The evaluation fibration $\omega_f \colon \map(X, Y;f)
\to Y$ is fibre-homotopy equivalent to $\omega_0 \colon \map(X, Y;0)
\to Y.$
\end{itemize}
\end{theorem}
\begin{proof}
The first statement follows from (\ref{eq:orbit=cyclic}) and
\thmref{thm:actions on Y}~(A). The  equivalence of (A) and (B) is a
consequence of adjointness, as mentioned above. We obtain (A)
implies (C) by observing that the universal action $A_\infty$ is
strictly unital and applying \thmref{thm:actions on Y}~(C) and
(\ref{eq:orbit=cyclic}). Finally, note that (C) implies (B) since
the evaluation fibration $\omega_0 \colon \map(X, Y;0) \to Y$ admits
the section $s(y)(x) = y.$
\end{proof}

\begin{corollary}
Let $Y$ be a finite CW complex.  Then $Y$ is an H-space if and only
if for every finite CW complex $X$ the evaluation fibrations
$\omega_f \colon \map(X, Y;f) \to Y$ are fibre-homotopy equivalent
for all maps $f \colon X \to Y.$
\end{corollary}
\begin{proof}
The result follows from \thmref{thm:cyclic} and the equivalences:
$$  \text{ $Y$ is  an H-space} \Leftrightarrow \text{ $1 \colon Y \to Y$ is cyclic}
\Leftrightarrow \text{every map $f \colon X \to Y$ is cyclic}$$
which are direct from  definitions.
\end{proof}

We now consider the above action $\mathcal{A}_\infty$ in the special
case in which $X$ is a co-H-space. Suppose the coproduct is $\sigma
\colon X \to X \vee X.$ The map $\sigma$ induces a pairing which we
denote `$+$' in the set $[X,Y]$. By \cite[Th.1.5]{Var} the set of
cyclic maps $G(X, Y)$ is a subgroup of $[X, Y]$  when $X$ is a
co-group-like space. When $X$ is merely a co-H-space, Varadarajan's
proof gives that the set $G(X, Y)$ is closed under addition. We show
that, when $X$ is a co-H-space,   the orbit of a class $[f] \in [X
,Y]$ under the action of $(\mathcal{A}_\infty)_\sharp$ is just the
set of translates of $[f]$ by $G(X, Y)$, that is, we have
\begin{equation}\label{eq:translates}
O_\infty([f]) = \{ d + f | d \in G(X, Y) \}
\end{equation}
This result is a direct consequence of the following:
\begin{lemma} \label{lemma:sum}
Let $X$ be a co-H-space.   Let $\gamma \colon X \to \map(Y, Y;1)$ and $f \colon X \to Y$
be based maps. Let $d \colon   X \to Y$ be defined by $d(x) = \gamma(x)(y).$
Then
$$\gamma \cdot f \sim_* d + f$$
\end{lemma}
\begin{proof}
Let $\Gamma \colon X \times Y \to Y$ denote the adjoint of $\gamma.$
By the definition of $d$ we
then have the following
homotopy-commutative diagram:
$$\xymatrix{X \ar[r]^-{\Delta} \ar[rd]_-{\sigma} & X \times X
\ar[r]^-{1\times f} & X \times Y \ar[r]^-{\Gamma} & Y\\
& X \vee X \ar[u]_{J} \ar[r]_-{1 \vee f} & X \vee Y \ar[u]^{J}
\ar[ru]_-{(d \mid 1)}& & \qed }$$
\renewcommand{\qed}{}
\end{proof}

The following consequence was proved by Yoon (\cite[Th.4.9]{Yoon89}) for $X$  a suspension.
\begin{theorem}\label{thm:suspension}
Suppose $X$ is a  CW co-H-space and $Y$ is any CW complex. Let $d \in
G(X,Y)$ be any cyclic map.  Then for each map $f
\colon X \to Y$, we have $\map(X,Y;f) \simeq \map(X,Y; f + d)$. If $X$
is a finite co-H-space then the corresponding evaluation
fibrations $\omega_f$ and $\omega_{f +d}$ are fibre homotopy equivalent.
\end{theorem}
\begin{proof}
Since $d \colon X \to Y$ is cyclic there exists a based map $\Gamma \colon
X \times Y \to Y$ extending $(d \mid 1) \colon X \vee Y \to Y$.
Let $\gamma \colon  X \to \map(Y, Y; 1)$ denote the adjoint to $\Gamma.$
Then, under the universal action we have $\gamma \circ f \sim_* d + f$
by \lemref{lemma:sum} and the result follows from
\thmref{thm:actions on Y}.
\end{proof}

Note that if $X$ is a suspension, or more generally a cogroup-like
space, then (\ref{eq:translates}) gives a bijection
\begin{equation}
\label{eq:coset} \mathcal{O}_\infty \cong [X, Y]/G(X, Y),
\end{equation}
where $\mathcal{O}_\infty$ denotes the orbits of the action
$(\mathcal{A}_\infty)_\sharp$ on $[X, Y]$ and the right-hand side is
simply the quotient group. If, for instance, $X = S^n$, then from
\corref{cor:surjection}~(C) we obtain a surjection
\begin{equation}\label{eq:Gottlieb surjection}
 \xymatrix{ \pi_n(Y)/ G_n(Y) \ar@{->>}[r] &
\begin{displaystyle}
\frac{\{ \text{evaluation fibrations   $\omega_f \colon
\map(S^n,Y;f) \to Y$} \}}{\text{fibre-homotopy equivalence}}
\end{displaystyle}.}
\end{equation}

When $Y$ is simple, the connecting homomorphism in the long exact
homotopy sequence of the evaluation fibration $\omega_f \colon
\map(S^n,Y;f) \to Y$, when viewed as a map
$$\partial \colon \pi_{k}(Y)  \to \pi_{k-1}(\map_*(S^n, Y;f) \cong
\pi_{k+n-1}(Y),$$
may be described in terms of Whitehead products with the class
represented by $f \colon S^n \to Y$ (see \cite[\S 3]{Wh}).  This
fact can be used in special cases to distinguish non-equivalent
components of $\map(S^n, Y)$ by calculating homotopy groups. For a
recent application of this method see \cite{Suth2}. Hansen uses this
method in \cite{Han1};    his result \cite[Th.2.3]{Han1} implies the
surjection (\ref{eq:Gottlieb surjection}) is actually a bijection
when $Y=S^m$ is also a sphere. Thus the homotopy classification of
components of $\map(S^n, S^m)$ reduces to the problem of computing
the Gottlieb groups $G_n(S^m)$. See \cite{G-M} for recent results in
this direction. Using other recent calculations of Gottlieb groups
we obtain the following.

\begin{example}\label{ex:Stiefel manifold}
Let $V_{n+k,k} = O(n+k)/O(n)$ denote the real Stiefel manifold.
Suppose that we have $m \leq n$ and $m \equiv 3,5,6,7 (\mod 8).$
Then the evaluation fibrations $$\omega_f \colon \map(S^m,
V_{n+k,k};f) \to V_{n+k,k}$$ \emph{are fibre-homotopy equivalent for
all $f \colon S^m \to V_{n +k, k}$.} For by \cite[Th.3.1]{LMW}, we
have $G_m(V_{n+k,k}) =\pi_m(V_{n+k,k})$ in these cases.  Similar
examples can be formulated for the complex and quaternionic Stiefel
manifolds using \cite[Th.3.2]{LMW} and \cite[Th.3.3]{LMW},
respectively.
\end{example}

Our last remark on these topics concerns the case in which $Y$ is a
so-called $G$-space, that is, a space that satisfies $G_n(Y) =
\pi_n(Y)$ for each $n$.  Such spaces have been studied by Siegel,
Gottlieb, and others, and are considered as being "close" to
$H$-spaces from certain points of view.  There are examples of
$G$-spaces that are not $H$-spaces, however \cite{Sie}. If $Y$ is a
$G$-space, then the surjection (\ref{eq:Gottlieb surjection}) yields
that all evaluation fibrations $\map(S^n,Y;f) \to Y$ are
fibre-homotopy equivalent to each other (and each has a section).
This is a further property that $G$-spaces share with $H$-spaces.

In case the function space of based maps is of interest, there is a
separate source of actions in addition to those obtained by
restricting actions on the unbased function space, as we have done
above. Namely, those group-like actions on the based mapping space
$\map_*(X, Y)$ that arise from cogroup-like actions on $X$.  We
finish the paper with a brief discussion of this topic. The approach
here is essentially that observed by Sutherland in \cite[\S
4]{Suth2}.

Suppose that $C$ is a \emph{co-H-space} with comultiplication
$\sigma \colon C \to C \vee C$.  Then, for any space $Y$, the based
function space $\map_*(C, Y)$ is an $H$-space with product $m(f, g)
= (f \mid g) \circ \sigma$, where $f, g \colon C \to Y$ are based
maps. By \cite[Th.III.5.16]{Wh}, $\map_*(C, Y)$ is
homotopy-associative, respectively group-like,  if $C$ is
homotopy-coassociative, respectively cogroup-like. By a
\emph{homotopy-coassociative coaction} of a homotopy-coassociative
co-H-space $C$ on a based space $X$, we mean a based map $B \colon X
\to C \vee X$ that satisfies $p_2 \circ B \sim_* 1\colon X \to X$
and $(\sigma \mid 1 ) \circ B \sim_* ( 1\mid B) \circ B  \colon X
\to C \vee C \vee X$ where $p_2\colon C \vee B \to B$ is the obvious
projection.

Suppose $C$ is a homotopy-coassociative co-H-space and $B \colon X
\to C \vee X$ is a homotopy-coassociative coaction.  Define
\begin{equation}\label{eq:action from coaction}
\mathcal{B} \colon \map_*(C, Y) \times \map_*(X, Y) \to \map_*(X,
Y).
\end{equation}
by setting $\mathcal{B}(\gamma, f) = (\gamma | f) \circ B.$ It is
direct to check that $\mathcal{B}$ defines a homotopy-associative
action on $\map_*(X,Y)$.  If $C$ is cogroup-like, this is a
group-like action. Thus we may apply \thmref{thm:components of an
action} to this situation, giving:

\begin{theorem}\label{thm:general result-coactions}
Let $B \colon X \to C \vee X$ be a homotopy-coassociative coaction
of a cogroup-like space $C$ on $X$. Let $\gamma \colon C \to Y$ be a
based map.
\begin{itemize}
\item[(A)] For any map $f \colon X \to Y$, under the resulting
action (\ref{eq:action from coaction}) on $\map_*(X,Y)$,
the path
components $\map_*(X,Y;f)$ and $\map_*(X,Y;\gamma\cdot f)$ of
$\map_*(X,Y)$ have the same homotopy type.
\item[(B)] Write
$$
\mathcal{B}_\# \colon [C, Y] \times [X, Y] \to [X, Y]
$$
for the action induced on homotopy sets by the action
(\ref{eq:action from coaction}). Let $\mathcal{O}'$ denote the set
of orbits of $[X,Y]$ under this action of the group $[C, X]$.  There
is a surjection of sets
$$\xymatrix{\mathcal{O}' \ar@{->>}[r] &
\begin{displaystyle}
\frac{\{ \text{components of $\map_*(X,Y)$} \}}{\simeq}
\end{displaystyle} .}$$
In particular, if $\mathcal{O}'$ is a finite set, then there are
finitely many distinct homotopy types amongst the components of
$\map_*(X,Y)$.
\end{itemize}
\end{theorem}

\begin{proof}
The proof is a direct consequence of the preceding discussion and
\thmref{thm:components of an action}.
\end{proof}

A standard source for a coaction on a space $X$ is a cofibration
sequence
$$\xymatrix{Z \ar[r]^{i} & A \ar[r]^{q} & X}$$
in which $X$ occurs as the cofibre.  For then we have $B\colon X \to
\Sigma Z \vee X$, the usual coaction of the cogroup-like space
$\Sigma Z$ on $X$. This then leads to an action as above
\begin{equation}\label{eq:action-specific cofibration}
\mathcal{B} \colon \map_*(\Sigma Z, Y) \times \map_*(X, Y) \to
\map_*(X, Y).
\end{equation}
Consider the Puppe sequence
\begin{equation}\label{eq:Puppe sequence cofibration}
\xymatrix{ \cdots \ar[r] & [\Sigma Z, Y] \ar[r]^{\partial^*} & [X,
Y] \ar[r]^{q^*} & [A, Y] \ar[r]^{i^*} & [Z, Y] }
\end{equation}
As is well-known, $[\Sigma Z,Y]$ acts on $[X, Y]$ and an orbit of $f
\in [X,Y]$ under this action is precisely the pre-image of $q^*(f)$
(see \cite[III.6.20]{Wh}).  It is easy to see that this action in
the Puppe sequence is identical with the action $\mathcal{B}_\#$
induced on $[X,Y]$ by the action (\ref{eq:action-specific
cofibration}). Thus the following result is a direct consequence of
\thmref{thm:general result-coactions}.

\begin{theorem}\label{thm:pre-image of q^*}
Let $q \colon A \to X$ be the cofibre projection of a cofibration in
which $X$ occurs as the cofibre.  Suppose $f\circ q \sim_* g\circ q
\colon A \to Y$ for maps $f, g \colon X \to Y$. Then $\map_*(X,Y;f)$
and $\map_*(X, Y; g)$ have the same homotopy type. If the image of
$q^* \colon [X,Y] \to [A, Y]$ is a finite set in $[A, Y]$, then
there are finitely many distinct homotopy types amongst the
components of $\map_*(X,Y)$. \qed
\end{theorem}

\begin{example}
Suppose $X$ is an $n$-dimensional manifold. Then $X$ occurs as the
cofibre in a cofibration of the form $S^n \to A \to X$ where $A$ is
an $(n-1)$-dimensional CW complex. Note  that the components of
$\map_*(X,S^n)$ are in one-to-one correspondence with $H^n(X)$, by
the Hopf-Whitney classification theorem, and so there are generally
infinitely many components of $\map_*(X, S^n)$.  However, $[A, S^n]$
consists of a single element, namely the homotopy class of the
trivial map.  By \thmref{thm:pre-image of q^*},    \emph{ all
components of the based mapping space $\map_*(X, S^n)$ have the same
homotopy type}.
\end{example}

Along the same lines, we offer the following:

\begin{example}
Suppose $X = S^n \cup_\alpha e^{r+1}$ is a two-cell complex, with
$\alpha \in \pi_r(S^n)$ for some $r\geq n$.  Suppose $Y$ is any
space with $\pi_n(Y)$ finite.  Then there are finitely many distinct
homotopy types amongst the components of the based mapping space
$\map_*(X, Y)$.  For we have a cofibre sequence $S^r \to S^n \to X$,
in whose Puppe sequence the map $q^* \colon [X, Y] \to [S^n, Y]$ has
finite image by hypothesis.  The assertion follows from
\thmref{thm:pre-image of q^*}.
\end{example}

\providecommand{\bysame}{\leavevmode\hbox
to3em{\hrulefill}\thinspace}
\providecommand{\MR}{\relax\ifhmode\unskip\space\fi MR }
\providecommand{\MRhref}[2]{%
  \href{http://www.ams.org/mathscinet-getitem?mr=#1}{#2}
} \providecommand{\href}[2]{#2}

\end{document}